\documentclass{article}

\usepackage{dsfont}
\usepackage[T1]{fontenc}
\usepackage[latin1]{inputenc}
\usepackage{amsmath}
\usepackage{amssymb}
\usepackage{amsthm}
\usepackage{amsfonts}
\usepackage{enumerate}
\usepackage{mathrsfs}
\usepackage[arrow, matrix, curve]{xy}
\usepackage{wrapfig}
\usepackage{graphicx}
\usepackage{pstricks,pst-node,pst-plot} %This loads the pstricks package
\usepackage{color} %This loads the color package
\usepackage{float}
\usepackage{mathdots}
\usepackage{nicefrac}
\usepackage{pstricks-add}
\usepackage{array,blkarray}
\usepackage[labelsep=none]{caption}

\newcommand{\edim}[0]{\operatorname{edim}}

\newcommand{\Ap}[0]{\operatorname{Ap}}

\makeatletter
\providecommand*{\cupdot}{%
  \mathbin{%
    \mathpalette\@cupdot{}%
  }%
}

\providecommand*{\bigcupdot}{%
  \mathop{%
    \vphantom{\bigcup}%
    \mathpalette\@bigcupdot{}%
  }%
}
\newcommand*{\@bigcupdot}[2]{%
  \ooalign{%
    $\m@th#1\bigcup$\cr
    \sbox0{$#1\bigcup$}%
    \dimen@=\ht0 %
    \advance\dimen@ by -\dp0 %
    \sbox0{\scalebox{2}{$\m@th#1\cdot$}}%
    \advance\dimen@ by -\ht0 %
    \dimen@=.5\dimen@
    \hidewidth\raise\dimen@\box0\hidewidth
  }%
}
\makeatother

\swapnumbers
\newtheorem{theorem}{Theorem}[section]

\newtheorem{proposition}[theorem]{Proposition}

\newtheorem{coro}[theorem]{Corollary}

\newtheorem*{exampleo}{Example}

\newtheorem*{definitiono}{Definition}

\newtheorem{remark}[theorem]{Remark}
\newtheorem{remarkTrick}[theorem]{}

\newtheorem*{remarko}{Remark}

\newcommand*{\hfillplus}{\hfill\linebreak[3]\hspace*{\fill}}

\author{M.~Hellus, A.~Rechenauer and R.~Waldi}
\title{Variants on a question of Wilf}
\date{\today}

\begin{document}

\maketitle

\begin{abstract}

Let $S\neq\mathds N$ be a numerical semigroup generated by $e$ elements. In his paper (A Circle-Of-Lights Algorithm for the ``Money-Changing Problem'', Amer. Math. Monthly 85 (1978), 562--565), H.~S.~Wilf raised the following question: Let $\Omega$ be the number of positive integers not contained in $S$ and $c-1$ the largest such element. Is it true that the fraction $\frac\Omega c$ of omitted numbers is at most $1-\frac1e$?

Let $B\subseteq\mathds N^{e-1}$ be the complement of an artinian $\mathds N^{e-1}$-ideal. Following a concept of A.~Zhai (An asymptotic result concerning a question of Wilf, arXiv:1111.2779v1 [math.CO]) we relate Wilf's problem to a more general question about the weight distribution on $B$ with respect to a positive weight vector. An affirmative answer is given in special cases, similar to those considered by R.~Fr\"oberg, C.~Gottlieb, R.~H\"aggkvist (On numerical semigroups, Semigroup Forum, Vol.~35, Issue 1, 1986/1987, 63--83) for Wilf's question.

\end{abstract}

\section{Averaging the weight of the points outside of an $\mathds N^{e-1}$-ideal}

\label{sect_aver}Let $e\in\mathds N_{\geq2}$ and $C\subseteq\mathds N^{e-1}$ be an $\mathds N^{e-1}$-\textbf{ideal}, i.\,e. $C+\mathds N^{e-1}\subseteq C$. Suppose that the corresponding monomial ideal $I(C)=(\{\underline X^c:=X_1^{c_1}\cdot\ldots\cdot X_{e-1}^{c_{e-1}}\vert c\in C\})\subseteq \mathds C[X_1,\ldots,X_{e-1}]$ is artinian. Then the complement $B=\mathds N^{e-1}\setminus C$ of $C$ is finite and $\{\underline X^b\vert b\in B\}$ is a vector space basis of the residue class ring $R(B):=\mathds C[X_1,\ldots, X_{e-1}]/I(C)$. Basic facts on monomial ideals can be found in \cite{HH}.

Choosing a \textbf{weight vector} $g=(g_1,\ldots,g_{e-1})\in\mathds R_{>0}^{e-1}$, the \textbf{weight} of the point $z=(z_1,\ldots,z_{e-1})\in\mathds Z^{e-1}$ is defined as the dot product $z\cdot g=z_1g_1+\ldots+z_{e-1}g_{e-1}$ of $z$ with $g$.

In \cite[Lemma 3]{Z}, Zhai has shown, that the mean weight of the elements of $B$ is bounded above by $\frac {e-1}e$ times their maximum weight, that is
\begin{equation}\label{Bound_Zhai}\frac1{\#B}\cdot\sum_{b\in B}b\cdot g\leq\frac {e-1}e\cdot\max(B\cdot g).\end{equation}
His proof even shows, that for the symmetric ($e-1$)-simplex
\[\Delta_{n,e}:=\{(x_1,\ldots,x_{e-1})\in\mathds N^{e-1}\vert x_1+\ldots+x_{e-1}\leq n-1\},\]
\begin{equation}\label{Simplex}\sum_{b\in\Delta_{n,e}}b\cdot g=\frac1e\#\,\Delta_{n,e}\cdot(n-1)(g_1+\ldots+g_{e-1}).\end{equation}

We shall consider $\mathds Z^{e-1}$ as a poset with regard to the canonical order
\[x\geq y\text{ if and only if }x-y\in\mathds N^{e-1}.\]
For $m\in\mathds N^{e-1}$ let $Q_m$ be the \textbf{cuboid}
\[Q_m:=\{b\in\mathds N^{e-1}\vert b\leq m\}.\]
Such cuboids are complements of ideals as well, and $m\cdot g=\max(Q_m\cdot g)$.

Let $m_1,\ldots,m_t$ be the maximum elements of $B$, hence $\bigcup_jQ_{m_j}=B$. Since $\{\underline X^{m_1},\ldots,\underline X^{m_t}\}$ induces a basis of the socle of the local ring $R(B)$, $t=t(R(B))$ is the \textbf{Cohen-Macaulay type} of $R(B)$. Similar to (\ref{Bound_Zhai}) we get

\begin{proposition}

\begin{enumerate}

\item[a)]
\begin{equation}\label{einseinsb}\frac1{\#B}\sum_{b\in B}b\cdot g\leq\frac{t(R(B))}{t(R(B))+1}\cdot \max(B\cdot g).\end{equation}

\item[b)] If $B=Q_m$ is a cuboid, then
\begin{equation}\label{Bound_Cube}\frac1{\#B}\cdot\sum_{b\in B}b\cdot g=\frac12\max(B\cdot g).\end{equation}

\end{enumerate}

\end{proposition}

\textbf{Proof} b) In fact $\frac12m$ is the center of symmetry of the cuboid $B=Q_m$, hence $\left(\frac12m\right)\cdot g$ is the mean weight of its lattice points.

a) $B=\bigcup_jQ_{m_j}$ and $Q_{m_j}=m_j-Q_{m_j}$ imply
\[\sum_{b\in B}b\cdot g\leq\sum_{j=1}^t\sum_{q\in Q_{m_j}}(m_j-q)\cdot g\leq t(\max(B\cdot g)\cdot\#B-\sum_{b\in B}b\cdot g).\]
\hfillplus$\square$

Obviously $B=Q_m$ if and only if $I(C)=(X_1^{\mu_1+1},\ldots,X_{e-1}^{\mu_{e-1}+1})$, $m=(\mu_1,\ldots,\mu_{e-1})$. Hence by \cite[Proposition A.6.5 and Corollary 1.3.6]{HH} we have

\begin{proposition} \label{prop_Gor}The following conditions are equivalent:

\begin{enumerate}

\item[a)] $B$ is a cuboid.

\item[b)] $I(C)$ is generated by pure powers of the variables $X_1,\ldots,X_{e-1}$.

\item[c)] $I(C)$ is a complete intersection ideal.

\item[d)] $R(B)$ is a Gorenstein ring.

\end{enumerate}

\end{proposition}\hfillplus$\square$

In section \ref{Apery} this will be applied to numerical semigroups.

\begin{remarkTrick}\textbf{A geometric interpretation of formula (\ref{Bound_Zhai}) in the sense of integral calculus}

Let $\beta=\max(B\cdot g)$ and $H_\beta\subseteq\mathds R^{e-1}$ be the hyperplane with the equation $g_1X_1+\ldots+g_{e-1}X_{e-1}=\beta$. Hence $B$ is contained in the half space $H_{\leq\beta}:g_1X_1+\ldots+g_{e-1}X_{e-1}\leq\beta$.

For $h\in\mathds R_{>0}$ let $\underline H\subseteq\mathds R^{e-1}\times\mathds R$ be the affine hull of $\{(0,h)\}\cup H_\beta\times\{0\}$. Above any point $b\in B$ there is exactly one point $(b,h_b)$ of $\underline H$ and $h_b\geq0$ is the ``height'' of $\underline H$ at the base point $b$.

In analogy to the formula for the volume of a pyramid we have
\begin{equation}\label{pyramide}\sum_{b\in B}h_b\geq \frac1e\cdot\#B\cdot h\ \text{(}=\frac1{\text{dimension}}\cdot\text{basis}\cdot\text{height}\text{)},\end{equation}
see figure 1, where the dotted lattice path denotes the boundary and the circles $\circ$ the maximum points of $B$.\hfillplus$\square$\end{remarkTrick}

(For (\ref{pyramide}) to be true it is essential, that $B$ is the complement of an ideal, i.\,e. that $b\in B$ implies $Q_b\subseteq B$.)

\textbf{Proof} of (\ref{pyramide}): Immediate from (\ref{Bound_Zhai}).\hfillplus$\square$

\begin{figure}[H]\centering
\begin{pspicture}(-.5,-3)(10,10)
\psline{->}(0,0)(8,4)
\psline{->}(0,0)(6,-2)
\psline{->}(0,0)(0,8)
\psline[linewidth=.5pt](2.5,0)(2.5,3.5)
\pspolygon[fillstyle=solid,fillcolor=white](0,7)(6,3)(4.5,-1.5)
\pspolygon[fillstyle=solid,fillcolor=gray](0,0)(1.5,-.5)(2.5,0)(1,.5)
\psdot(2.5,3.5)
\uput[225](0,0){$0$}
\uput[180](0,7){$(0,h)$}
\uput[0](2.5,0){$b$}
\uput[0](.35,-.5){$Q_b$}
\uput[0](2.5,3.5){$(b,h_b)$}
\uput[0](3.25,3){$\underline H$}
\uput[0](5.25,.75){$H_\beta$}
\uput[0](0,7.5){$\mathds R$}
\uput[0](7,1){$\mathds R^{e-1}$}
\uput[0](3.1,-.2){$B$}
\psline[linestyle=dotted](4,2)(3.5,1.75)(5,1.25)(3.5,.5)(5,0)(3,-1)
\psdot[dotstyle=o,dotsize=5pt 0](4,2)
\psdot[dotstyle=o,dotsize=5pt 0](5,1.25)
\psdot[dotstyle=o,dotsize=5pt 0](5,0)
\end{pspicture}
\caption{}
\end{figure}

\section{Ap\'{e}ry sets}
\label{Apery}

\subsection{Zhai's version of Wilf's inequality}

Let $S=\mathds N\cdot g_0+\ldots+\mathds N\cdot g_{e-1}\neq\mathds N$ be a numerical semigroup with minimal generating set $\{g_0,\ldots,g_{e-1}\}\subseteq\mathds N_{>0}$,
\[g_0<\ldots<g_{e-1}, \gcd(g_0,\ldots,g_{e-1})=1.\]
We set
\[c=c(S):=\min\{s\in S\vert s+\mathds N\subseteq S\},\]
sometimes called the \textbf{conductor} of $S$ and
\[L=L(S):=S\cap\{0,\ldots,c-1\}\]
the part of $S$ to the \textbf{left} of $c$.
\[A=\Ap(S,g_0):=\{s\in S\vert s-g_0\not\in S\}\]
is called the \textbf{Ap\'{e}ry set} of $S$ with respect to $g_0$.

Each $a\in A$ can be written in the form $a=x\cdot g$, $g:=(g_1,\ldots,g_{e-1})$, $x\in\mathds N^{e-1}$.

Then each $s\in S$ has a unique presentation
\begin{equation}\label{uniq_pres}s=a_0g_0+a\text{, }a_0\in\mathds N\text{ and }a\in A.\end{equation}
In the following let us consider the \textbf{question of Wilf} from \cite{W}, which asks if the inequality
\begin{equation}\label{Wilf}\frac{\#L}c\geq\frac1e\end{equation}
holds for every numerical semigroup $S$.

With the help of his lemma, loc.\,cit., A.~Zhai succeeded in proving a weakened version of formula (\ref{Wilf}). Here we shall repeat the results and arguments of Zhai, as far as they seem to be useful in our later considerations.

Following Zhai we endow $\mathds N^{e-1}$ with the \textbf{(purely) lexicographic order} LEX, i.\,e. $a<_\text{LEX}b$ if the leftmost nonzero component of $a-b$ is negative.

For $a\in A$ let $\tilde a\in\mathds N^{e-1}$ be the LEX-\textbf{minimal} element
\[\tilde a:=\min_{\text{LEX}}\{x\in\mathds N^{e-1}\vert x\cdot g=a\}\]
of $\mathds N^{e-1}$ with weight $\tilde a\cdot g=a$. Set $\tilde A:=\{\tilde a\vert a\in A\}\subseteq\mathds N^{e-1}$.

Let
\[\pi:\mathds Z\times\tilde A\to\mathds Z, x\mapsto x\cdot (g_0,g)\]
be the restriction of the dot product $\mathds R^e\to\mathds R, x\mapsto x\cdot (g_0,g)$ to $\mathds Z\times\tilde A$. Because of the uniqueness of the presentation (\ref{uniq_pres}) and the construction of $\tilde A$,

\begin{itemize}

\item $\pi$ is bijective

\item $\pi(\mathds N\times\tilde A)=S$

\end{itemize}

According to Zhai $\tilde A\subseteq\mathds N^{e-1}$ is the complement of an ideal (\cite[Proof of theorem 1]{Z}). Since $\#\tilde A=g_0$, $\tilde a\cdot g=a$ for $a\in A$ and $c+g_0-1=\max A$, Zhai's inequality (\ref{Bound_Zhai}) implies
\begin{equation}\label{sechs_neu}\frac1{g_0}\sum_{a\in A}a\leq \frac{e-1}e(c+g_0-1).\end{equation}
In \cite[Lemma 1]{Z} one can find the exact formula
\begin{equation}\label{sechs}\#L=\frac1{g_0}\sum_{a\in A}(c-a)+\frac12(g_0-1)\end{equation}
for the cardinality of $L$. Hence Wilf's inequality (\ref{Wilf}) reads as
\begin{equation}\label{acht_neu}\frac1{g_0}\sum_{a\in A}a\leq\frac{e-1}ec+\frac12(g_0-1)\end{equation}
which for $e\geq3$ is stronger than Zhai's version (\ref{sechs_neu}).

One obtains

\begin{coro} The numerical semigroup $S$ is symmetric if and only if formula (\ref{Bound_Cube}) holds for $B=\tilde A$.\end{coro}

\textbf{Proof} Here we have to take $\#B=\#\tilde A=g_0$ and $(g_1,\ldots,g_{e-1})$ as weight vector $g$. Further $B\cdot g=A$, $\max A=c+g_0-1$ and $\tilde a\cdot g=a$ for $\tilde a\in\tilde A$. Hence formula (\ref{Bound_Cube}) means
\[\frac1{g_0}\sum_{a\in A}a=\frac12(c+g_0-1),\]
i.\,e.
\begin{equation}\label{acht}\frac12c=\frac1{g_0}\sum_{a\in A}(c-a)+\frac12(g_0-1).\end{equation}
But by (\ref{sechs}) the right hand side of (\ref{acht}) always equals $\#L$. Hence (\ref{Bound_Cube}) holds if and only if $\frac12c=\#L$, i.\,e. if $S$ is symmetric.\hfillplus$\square$

\subsection{Comparing the types of $S$ and of $R(\tilde A)$}

Let $\mathds C[S]=\mathds C[\{t^s|s\in S\}]\subseteq\mathds C[t]$  be the semigroup ring of $S$ and $R(\tilde A)=\mathds C[X_1,\ldots,X_{e-1}]/I(\mathds N^{e-1}\setminus\tilde A)$.

We denote by $t(R)$ the type of the local Cohen-Macaulay ring $(R,m)$. If $R$ is artinian then $t(R)$ can be computed as the vector space dimension of the socle $0:_Rm$ of $R$. Further the type $t(S):=t(\mathds C[[S]])$ of a numerical semigroup $S$ equals the number of pseudo-Frobenius numbers of $S$, these are the elements $f\in\mathds N\setminus S$ such that $f+S\setminus\{0\}\subseteq S$.

\begin{proposition}\label{Comp_types}
\[t(R(\tilde A))\geq t(S).\]
\end{proposition}

\textbf{Proof} The monomials $\underline X^s$, $s$ a maximum element of $\tilde A$, induce a basis of the socle of $R(\tilde A)$. Further the preimage $\pi^{-1}(f)$ of any pseudo-Frobenius number $f$ of $S$ is a maximum element of $\{-1\}\times\tilde A\subseteq\mathds Z^e$, hence $t(R(\tilde A))\geq t(S)$.\hfillplus$\square$

\begin{proposition}\label{Prop_in_ideal} $I(\mathds N^{e-1}\setminus\tilde A)$ is the initial ideal $\operatorname{In}_{\text{LEX}}(J)$ of
\[J=\ker(\mathds C[X_1,\ldots,X_{e-1}]\twoheadrightarrow\mathds C[S]/(t^{g_0}), X_i\mapsto t^{g_i}\mod(t^{g_0})),\]
with respect to LEX.\end{proposition}

\textbf{Proof} Let $I(S)$ be the kernel of
\[\mathds C[X_0,\ldots,X_{e-1}]\twoheadrightarrow\mathds C[S], X_i\mapsto t^{g_i}.\]
Since
\[S=\mathds N\cdot g_0+\tilde A\cdot g,\]
\[I(S)=(\{\underline X^\alpha-\underline X^\beta\vert\alpha\in\mathds N\times\tilde A, \beta\in\{0\}\times(\mathds N^{e-1}\setminus\tilde A)\text{ such that }\alpha\cdot(g_0,g)=\beta\cdot(g_0,g)\in S\}).\]
Hence
\[I(\mathds N^{e-1}\setminus\tilde A)=(\{X_1^{\beta_1}\cdot\ldots\cdot X_{e-1}^{\beta_{e-1}}\vert(\beta_1,\ldots,\beta_{e-1})\in(\mathds N^{e-1}\setminus\tilde A)\})\subseteq\operatorname{In}_{\text{LEX}}(J).\]
Since
\[\dim_{\mathds C}\mathds C[X_1,\ldots,X_{e-1}]/I(\mathds N^{e-1}\setminus\tilde A)=\#\tilde A=g_0=\]
\[=\dim_{\mathds C}\mathds C[X_1,\ldots,X_{e-1}]/J=\dim_{\mathds C}\mathds C[X_1,\ldots,X_{e-1}]/\operatorname{In}_{\text{LEX}}(J),\]
even $\operatorname{In}_{\text{LEX}}(J)=I(\mathds N^{e-1}\setminus\tilde A)$.\hfillplus$\square$

Using Proposition \ref{prop_Gor} (\cite[loc.\,cit.]{HH}) we get another, sufficient, condition on $\tilde A$ for $S$ being symmetric:

\begin{coro}

\begin{enumerate}

\item[a)] If $\tilde A$ is a cuboid, then $S$ is a complete intersection, in particular symmetric.

\item[b)] There are complete intersections in embedding dimension $3$, with $\tilde A$ not being a rectangle, in particular $t(R(\tilde A))>t(S)$.

\end{enumerate}

\end{coro}

\textbf{Proof}

\begin{enumerate}

\item[a)]
By \ref{prop_Gor} and \ref{Prop_in_ideal}, $\operatorname{In}_{\text{LEX}}(J)=I(\mathds N^{e-1}\setminus\tilde A)$ is a complete intersection, and so is $J$, and finally $\mathds C[S]$.

\item[b)] \textbf{Examples} (i) $S=\langle7,8,12\rangle$ is symmetric, hence a complete intersection, since $e=3$:
\[L=\{0,7,8,12,14,15,16,19,20,21,22,23,24\}\]
\[A=\{0,8,12,16,20,24,32\}\]
In the following pictures, we mark the elements of $\tilde A$ by $\bullet$:

\begin{figure}[H]\centering
\begin{pspicture}(-1,-1)(4.5,4)
\psline{->}(0,0)(3.5,0)
\psline{->}(0,0)(0,3.5)
\psdot[dotsize=5pt](0,0)
\uput[225](0,0){$0$}
\psdot[dotsize=5pt](1,0)
\uput[225](1,0){$8$}
\psdot[dotsize=5pt](2,0)
\uput[225](2,0){$16$}
\psdot[dotsize=5pt](0,1)
\uput[225](0,1){$12$}
\psdot[dotsize=5pt](1,1)
\uput[225](1,1){$20$}
\psdot[dotsize=5pt](0,2)
\uput[225](0,2){$24$}
\psdot[dotsize=5pt](1,2)
\uput[225](1,2){$32$}
\psdot[dotsize=2pt](3,0)
\psdot[dotsize=2pt](2,1)
\psdot[dotsize=2pt](3,1)
\psdot[dotsize=2pt](2,2)
\psdot[dotsize=2pt](3,2)
\psdot[dotsize=2pt](0,3)
\psdot[dotsize=2pt](1,3)
\psdot[dotsize=2pt](2,3)
\psdot[dotsize=2pt](3,3)
\end{pspicture}
\caption{: $\tilde A$ is not a rectangle}
\end{figure}
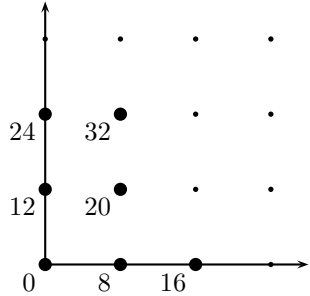

(ii) $S=\langle6,7,9\rangle$ is symmetric, hence a complete intersection, since $e=3$:
\[L=\{0,6,7,9,12,13,14,15,16\}\]
\[A=\{0,7,9,14,16,23\}\]

\begin{figure}[H]\centering
\begin{pspicture}(-1,-1)(4.5,3)
\psline{->}(0,0)(3.5,0)
\psline{->}(0,0)(0,2.5)
\psdot[dotsize=5pt](0,0)
\uput[225](0,0){$0$}
\psdot[dotsize=5pt](1,0)
\uput[225](1,0){$7$}
\psdot[dotsize=5pt](2,0)
\uput[225](2,0){$14$}
\psdot[dotsize=5pt](0,1)
\uput[225](0,1){$9$}
\psdot[dotsize=5pt](1,1)
\uput[225](1,1){$16$}
\psdot[dotsize=5pt](2,1)
\uput[225](2,1){$23$}
\psdot[dotsize=2pt](3,0)
\psdot[dotsize=2pt](3,1)
\psdot[dotsize=2pt](0,2)
\psdot[dotsize=2pt](1,2)
\psdot[dotsize=2pt](2,2)
\psdot[dotsize=2pt](3,2)
\end{pspicture}
\caption{: $\tilde A$ is a rectangle}
\end{figure}
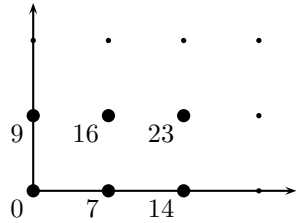

\hfillplus$\square$

\end{enumerate}

Next we will see how $\mathds Z^{e-1}$ can be tesselated by its subset $\tilde A$.

Let $B$ as in section \ref{sect_aver} be the complement of an artinian ideal $C\subseteq\mathds N^{e-1}$.

\begin{definitiono} $B$ induces a \textbf{periodic tesselation} on $\mathds Z^{e-1}$ if there exists an $(e-1)$-dimensional sublattice $\Lambda\subseteq\mathds Z^{e-1}$ such that
\[\mathds Z^{e-1}=\bigcupdot_{\lambda\in\Lambda}(\lambda+B).\]
\end{definitiono}

\begin{proposition}\label{tessi}$\tilde A$ induces a periodic tesselation on $\mathds Z^{e-1}$.\end{proposition}

In example (i), $S=\langle7,8,12\rangle$, we have
\[\Lambda=\mathds Z\cdot(2,1)\oplus\mathds Z\cdot(-3,2).\]

\begin{figure}[H]\centering
\begin{pspicture}(-3,-1.5)(3.5,3)
\psdot[dotsize=2pt](0,0)
\psdot[dotsize=2pt](.5,0)
\psdot[dotsize=2pt](1,0)
\psdot[dotsize=2pt](0,.5)
\psdot[dotsize=2pt](.5,.5)
\psdot[dotsize=2pt](0,1)
\psdot[dotsize=2pt](.5,1)
\psframe[fillcolor=lightgray,fillstyle=solid,linestyle=none,linewidth=0](-2.5,.5)(-2,1.5)
\psline(-2.5,.5)(-1.5,.5)
\psline(-2.5,.5)(-2.5,1.5)
\psline(-2.5,1.5)(-2,1.5)
\psline(-2,1.5)(-2,.5)
\psframe[fillcolor=lightgray,fillstyle=solid,linestyle=none,linewidth=0](-1.5,1)(-1,2)
\psline(-1.5,1)(-.5,1)
\psline(-1.5,1)(-1.5,2)
\psline(-1.5,2)(-1,2)
\psline(-1,2)(-1,1)
\psframe[fillcolor=lightgray,fillstyle=solid,linestyle=none,linewidth=0](-.5,1.5)(0,2.5)
\psline(-.5,1.5)(.5,1.5)
\psline(-.5,1.5)(-.5,2.5)
\psline(-.5,2.5)(0,2.5)
\psline(0,2.5)(0,1.5)
\psframe[fillcolor=lightgray,fillstyle=solid,linestyle=none,linewidth=0](1,.5)(1.5,1.5)
\psline(1,.5)(2,.5)
\psline(1,.5)(1,1.5)
\psline(1,1.5)(1.5,1.5)
\psline(1.5,1.5)(1.5,.5)
\psframe[fillcolor=lightgray,fillstyle=solid,linestyle=none,linewidth=0](2.5,-.5)(3,.5)
\psline(2.5,-.5)(3.5,-.5)
\psline(2.5,-.5)(2.5,.5)
\psline(2.5,.5)(3,.5)
\psline(3,.5)(3,-.5)
\psframe[fillcolor=lightgray,fillstyle=solid,linestyle=none,linewidth=0](1.5,-1)(2,0)
\psline(1.5,-1)(2.5,-1)
\psline(1.5,-1)(1.5,0)
\psline(1.5,0)(2,0)
\psline(2,0)(2,-1)
\psframe[fillcolor=lightgray,fillstyle=solid,linestyle=none,linewidth=0](.5,-1.5)(1,-.5)
\psline(.5,-1.5)(1.5,-1.5)
\psline(.5,-1.5)(.5,-.5)
\psline(.5,-.5)(1,-.5)
\psline(1,-.5)(1,-1.5)
\psframe[fillcolor=lightgray,fillstyle=solid,linestyle=none,linewidth=0](-1,-.5)(-.5,.5)
\psline(-1,-.5)(0,-.5)
\psline(-1,-.5)(-1,.5)
\psline(-1,.5)(-.5,.5)
\psline(-.5,.5)(-.5,-.5)
\pscircle(0,0){.1}
\pscircle(-2.5,.5){.1}
\pscircle(-1.5,1){.1}
\pscircle(-.5,1.5){.1}
\pscircle(1,.5){.1}
\pscircle(2.5,-.5){.1}
\pscircle(1.5,-1){.1}
\pscircle(.5,-1.5){.1}
\pscircle(-1,-.5){.1}
\psline[linestyle=dashed]{->}(0,0)(1,.5)
\psline[linestyle=dashed]{->}(0,0)(-1.5,1)
\end{pspicture}
\caption{: $\Lambda=\{\circ\}, \tilde A=\{\cdot\}$}
\end{figure}
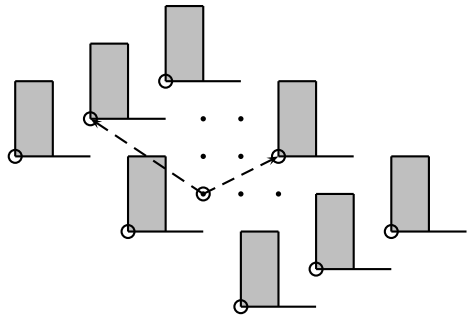

\textbf{Proof of \ref{tessi} Proposition} Let
\[\begin{array}{cccccc}\varepsilon:&\mathds Z^{e-1}&\buildrel\tau\over\to&\mathds Z&\buildrel\text{can.}\over\to&\mathds Z/g_0\mathds Z\\&x&\mapsto&x\cdot g&\mapsto&x\cdot g\mod g_0\end{array}\]
Then $\Lambda:=\ker\varepsilon\subseteq\mathds Z^{e-1}$ is an $(e-1)$-dimensional lattice. By the definition of $\tilde A$ and $A=\Ap(S,g_0)$ we have $A=\tau(\tilde A)$, and $\varepsilon$ maps $\tilde A\subseteq\mathds Z^{e-1}$ bijectively onto $\mathds Z/g_0Z$, hence
\[\mathds Z^{e-1}=\bigcupdot_{\tilde a\in\tilde A}(\tilde a+\Lambda)=\bigcupdot_{\lambda\in\Lambda}(\lambda+\tilde A).\]
\hfillplus$\square$

For $x=(x_1,\ldots,x_{e-1})\in\mathds R^{e-1}$ let $\operatorname{supp}(x):=\{i\vert i\in\{1,\ldots,e-1\},x_i\neq0\}$.

\begin{proposition}\label{prop_tess_equiv}Let $B$ be the complement of an artinian ideal $C\subseteq\mathds N^{e-1}$. Suppose $B$ induces a periodic tesselation $\mathds Z^{e-1}=\bigcupdot_{\lambda\in\Lambda}(\lambda+B)$.

\begin{enumerate}

\item[a)] If $x$ and $x'$ are minimal elements of $C$ such that $x-x'\in\Lambda$ and $\operatorname{supp}(x)\cap\operatorname{supp}(x')\neq\emptyset$, then $x=x'$.

\item[b)] There is at most one minimal element $x$ of $C$ outside the coordinate hyperplanes, and then in addition $x\in\Lambda$.

\end{enumerate}

\end{proposition}

\textbf{Proof} a) The canonical map $\varepsilon:\mathds Z^{e-1}\to\mathds Z^{e-1}/\Lambda$ operates bijectively on $B$, since $\mathds Z^{e-1}=\bigcupdot_{b\in B}(b+\Lambda)$. Let $i\in\operatorname{supp}(x)\cap\operatorname{supp}(x')$ and $e_i\in\mathds Z^{e-1}$ the $i$-th unit vector. Then $\{x-e_i,x'-e_i\}\subseteq B$ and $\varepsilon(x-e_i)-\varepsilon(x'-e_i)=\varepsilon(x-x')=0$. Hence $x-e_i=x'-e_i$, since $\varepsilon\vert_B$ is injective, and $x=x'$.

b) Let $x\in\mathds N_{>0}^{e-1}$ be a minimal element of $C$. By a) it suffices to show that $x\in\Lambda$: Since $\operatorname{supp}(x)=\{1,\ldots,e-1\}$,
\begin{equation}\label{zehn_neu}x-e_j\in B\text{ for }j=1,\ldots,e-1.\end{equation}
Since $\mathds Z^{e-1}=\bigcup_{\lambda\in\Lambda}(\lambda+B)$ and $x\not\in B$ there is an $\lambda\in\Lambda\setminus\{0\}$ such that $x-\lambda=:b\in B$. From $B\cap(\lambda+B)=\emptyset$ and (\ref{zehn_neu}) we conclude
\[b-e_j=(x-e_j)-\lambda\not\in B, j=1,\ldots,e-1,\]
hence $b=0$ and $x=\lambda\in\Lambda$; notice that $C$ is an ideal of $\mathds N^{e-1}$.\hfillplus$\square$

\begin{coro}\label{coro_tess_three_equiv}If $e=3$, then the following are equivalent:

\begin{enumerate}

\item[a)] $B$ induces a periodic tesselation on $\mathds Z^2$.

\item[b)] $t(R(B))\leq2$.

\end{enumerate}

\end{coro}

\textbf{Proof} a)$\Rightarrow$b): By \ref{prop_tess_equiv}, $C$ has at most one minimal element outside the axes of $\mathds R^2$; hence $B$ has at most two maximum elements, i.\,e. $t(R(B))\leq2$.

b)$\Rightarrow$a): If $t(R(B))=1$, then $B$ is a rectangle (\ref{prop_Gor} Proposition). In case $t(R(B))=2$ let $(x_1,y_1)$ and $(x_2,y_2)$ be the maximum elements of $B$, $x_1<x_2$ and $y_1>y_2$. Then $\Lambda=\mathds Z\cdot(x_1+1,y_2+1)\oplus\mathds Z\cdot(-x_2-1,y_1-y_2)$ works.\hfillplus$\square$

Applying \ref{tessi} and \ref{coro_tess_three_equiv} to numerical semigroups $S$ we obtain:

\begin{coro}\label{coro_33_neu}If $\edim S=3$, then $t(R(\tilde A))\leq2$. In particular, using \ref{Comp_types} Proposition, we obtain
\[t(S)=t(R(\tilde A))=2\text{ if }S\text{ is not symmetric.}\]
Hence $\tilde A$ is in the form of an ``$L$''.\end{coro}

See also \cite[theorem 11]{FGH} and its proof for similar considerations.

\begin{remarko} For $e\geq4$, $\tilde A$ can be of a more complicate ``staircase shape''.\end{remarko}

\begin{exampleo} Let $S=\langle9,10,12,13\rangle$. Both $S$ and $R(\tilde A)$ have type $5$:\end{exampleo}

\begin{figure}[H]\centering
\psset{unit=.3}
\begin{pspicture}(-6,-4)(7,9)
\psline(0,0)(0,6)
\psline(0,0)(4,0)
\psline(0,0)(-2,-2)
\psline(2,0)(2,4)
\psline(4,0)(4,2)
\psline(-1,-1)(-1,3)
\psline(-2,-2)(-2,0)
\psline(0,2)(4,2)
\psline(0,2)(-2,0)
\psline(0,4)(2,4)
\psline(0,4)(-1,3)
\psline(-2,-2)(-4,-2)
\psline(-2,0)(-4,0)
\psline(-1,1)(-3,1)
\psline(-1,3)(-3,3)
\psline(0,4)(-2,4)
\psline(0,6)(-2,6)
\psline(-2,6)(-2,4)
\psline(-2,4)(-3,3)
\psline(-3,3)(-3,1)
\psline(-3,1)(-4,0)
\psline(-4,0)(-4,-2)
\psline(4,0)(5,1)
\psline(5,1)(5,3)
\psline(5,3)(3,3)
\psline(3,3)(3,5)
\psline(3,5)(1,5)
\psline(1,5)(1,7)
\psline(1,7)(-1,7)
\psline(-1,7)(-2,6)
\psline(4,2)(5,3)
\psline(2,2)(3,3)
\psline(2,4)(3,5)
\psline(0,4)(1,5)
\psline(0,6)(1,7)
\psline{->}(-.5,6.5)(-.5,8)
\psline{->}(4.5,1.5)(6,1.5)
\psline{->}(-3,-1)(-4.5,-2.5)
\psdot[dotsize=2pt](-.5,6.5)
\psdot[dotsize=2pt](4.5,1.5)
\psdot[dotsize=2pt](-3,-1)
\uput[225](-4.5,-2.5){$x_1$}
\uput[0](6,1.5){$x_2$}
\uput[90](-.5,8){$x_3$}
\end{pspicture}
\caption{}
\label{fig_staircase}
\end{figure}

The figure indicates the $9$ unit cubes centered in the points of $\tilde A$. The $5$ steps correspond to the monomial basis of the socle of $R(\tilde A)$. The lattice $\Lambda$ is
\[\Lambda=\ker\varepsilon=\mathds Z\cdot(0,3,0)\oplus\mathds Z\cdot(3,-1,0)\oplus\mathds Z\cdot(1,1,-1).\]
More generally, for $S=S(n)=\langle n^2,n^2+1,n^2+n,n^2+n+1\rangle$, both $S$ and $R(\tilde A)$ have type $2n-1$ and formula (\ref{Wilf}) holds (cf. \cite{HRW}). Like in figure \ref{fig_staircase}, $\tilde A$ is a double staircase for $n$ stories, i.\,e.
\[\tilde A=B_n:=\{(x_1,x_2,x_3)\in\mathds N^3\vert x_1+x_2+x_3\leq n-1, x_1x_2=0\}.\]
We will even show

\begin{remark}

For all numerical semigroups $S=\langle g_0,g_1,g_2,g_3\rangle$ with $\tilde A=B_n$, $n\geq2$, Wilf's question has a positive answer.

\end{remark}

\textbf{Proof} Here $g_0=\#\,B_n=n^2$ and $c+g_0-1=\max A=(n-1)g_3$. By \cite{E}, formula (\ref{Wilf}) is true if $c\leq3g_0$. Hence it suffices to show, that (\ref{acht_neu}) holds in case $(n-1)g_3\geq4n^2$. For $B_{n,i}:=B_n\cap\{x_i=0\}$, $i=1,2$
\[B_{n,1}\cup B_{n,2}=B_n\text{ and }B_{n,1}\cap B_{n,2}=\{0,1,\ldots,n-1\}e_3.\]

Applying formula (\ref{Simplex}) to the triangles $B_{n,1}$ and $B_{n,2}$ yields inequality (\ref{acht_neu}) for $S$, if $(n-1)g_3\geq4n^2$:
\begin{align*}\sum_{B_n}b\cdot g&\leq\sum_{B_n}b\cdot(g_3-2,g_3-1,g_3)\\&=\sum_{i=1}^2\sum_{B_{n,i}}b\cdot(g_3-2,g_3-1,g_3)-\frac12n(n-1)g_3\\&=\frac23(n^2+n)(n-1)g_3-\frac12(n^2+n)(n-1)-\frac12n(n-1)g_3\\&=\frac34n^2(n-1)g_3-\frac1{12}(n-2)n(n-1)g_3-\frac12n(n^2-1)\\&\leq\frac34n^2(n-1)g_3-\frac13n^2(n-2)n-\frac12n(n^2-1)\text{, since }(n-1)g_3\geq4n^2\text{,}\\&\leq\frac34n^2(n-1)g_3-\frac14n^2(n^2-1)\text{, since }n\geq 2\text{,}\\&=n^2\left(\frac34c+\frac12(n^2-1)\right).\end{align*}
\hfillplus$\square$

\begin{exampleo} For $S=\langle9,20,21,23\rangle$ we have $\tilde A=B_3$, whereas the hypothesis of \cite[2.2\,Cor.]{HRW} is not fulfilled.\end{exampleo}

\begin{remarkTrick}\textbf{The Wilf ratio of a semigroup}

Recall that $\pi:\mathds Z\times\tilde A\to\mathds Z, x\mapsto x\cdot (g_0,g)$ is bijective and maps $\mathds N\times\tilde A$ onto $S$.

For $a\in\mathds Z$ let $H_{\geq a}:=\pi^{-1}(\mathds Z_{\geq a})$, similarly $H_{\leq a}$ and $H_a$. Hence under the bijection $\pi$ the following sets correspond to each other:
\[H_{\geq0}\widehat{=}\mathds N\]
\[H_{\leq c-1}\widehat{=}\mathds Z_{\leq c-1}\]
\[H:=H_{\geq0}\cap H_{\leq c-1}\widehat{=}\{0,\ldots,c-1\}\]
\[U:=H\cap\{x_0\geq0\}\widehat{=}L\]
\[D:=H\cap\{x_0\leq-1\}\widehat{=}\mathds N\setminus S\]
Hence Wilf's (proposed) inequality means that the \textbf{Wilf ratio of} $\mathbf S$
\[w(S):=\frac{\#U}{\#D}\]
is at least $\frac1{e-1}$.\hfillplus$\square$\end{remarkTrick}

\section{A more general question}

\subsection{Rephrasing the problem}

\label{subsect31}With the notation of section 2, Wilf asked if $(e-1)\#U\geq\#D$. Hence we will look at the subsets $U$ and $D$ of $H$ more precisely.
Let $n_0:=\lfloor\frac{c-1+g_0}{g_0}\rfloor=\lfloor\frac{\max A}{g_0}\rfloor$, the greatest integer such that $n_0g_0\leq c-1+g_0$.

For $i\in\mathds Z$ and $\tilde a\in\tilde A$,
\begin{equation}\label{neun}(i,\tilde a)\in H\text{ if and only if }0\leq ig_0+\tilde a\cdot g\leq c-1.\end{equation}

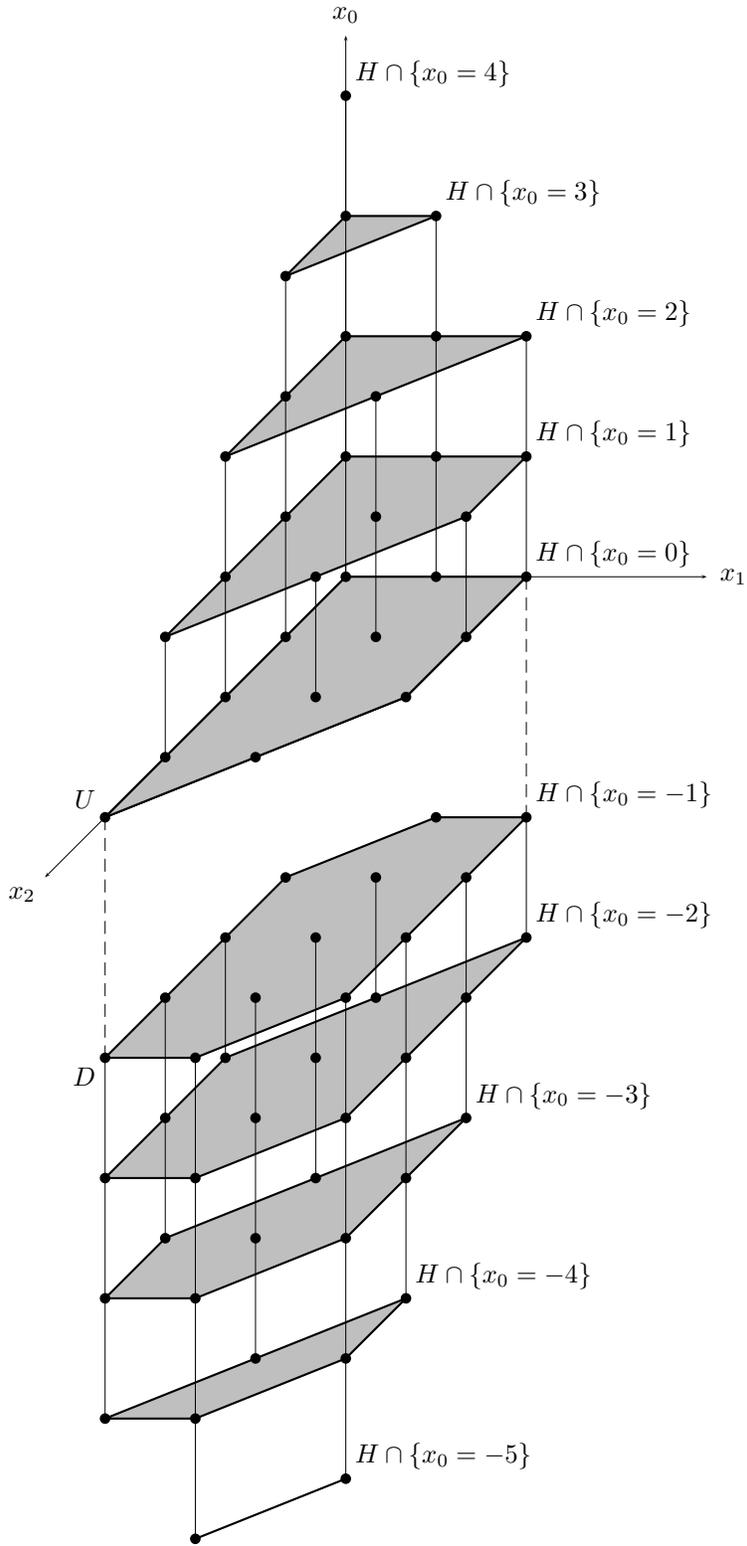
\begin{figure}[H]\centering
\psset{unit=.8}\begin{pspicture}(-7,-16.5)(8,7.5)
\psline[linewidth=.0001pt]{->}(0,0)(0,9)
\psline[linewidth=.0001pt]{->}(0,0)(6,0)
\psline[linewidth=.0001pt]{->}(0,0)(-5,-5)
\uput[90](0,9){$x_0$}
\uput[0](6,0){$x_1$}
\uput[225](-5,-5){$x_2$}
\pspolygon[fillstyle=solid,fillcolor=lightgray](0,6)(1.5,6)(-1,5)
\pspolygon[fillstyle=solid,fillcolor=lightgray](0,4)(3,4)(-2,2)
\pspolygon[fillstyle=solid,fillcolor=lightgray](-3,-1)(2,1)(3,2)(0,2)
\pspolygon[fillstyle=solid,fillcolor=lightgray](-4,-4)(1,-2)(3,0)(0,0)
\pspolygon[fillstyle=solid,fillcolor=lightgray](-4,-4)(1,-2)(3,0)(0,0)
\pspolygon[fillstyle=solid,fillcolor=lightgray](1.5,-4)(3,-4)(0,-7)(-2.5,-8)(-4,-8)(-1,-5)
\pspolygon[fillstyle=solid,fillcolor=lightgray](3,-6)(0,-9)(-2.5,-10)(-4,-10)(-2,-8)
\pspolygon[fillstyle=solid,fillcolor=lightgray](2,-9)(0,-11)(-2.5,-12)(-4,-12)(-3,-11)
\pspolygon[fillstyle=solid,fillcolor=lightgray](1,-12)(0,-13)(-2.5,-14)(-4,-14)
\psline(-2.5,-16)(0,-15)
\psline[linewidth=.0001pt](-3,-3)(-3,-1)
\psline[linewidth=.0001pt](-2,-2)(-2,2)
\psline[linewidth=.0001pt](-1,-1)(-1,5)
\psline[linewidth=.0001pt](0,0)(0,8)
\psline[linewidth=.0001pt](1.5,0)(1.5,6)
\psline[linewidth=.0001pt](3,0)(3,4)
\psline[linewidth=.0001pt](-.5,-2)(-.5,0)
\psline[linewidth=.0001pt](.5,-1)(.5,3)
\psline[linewidth=.0001pt](2,-1)(2,1)
\psline[linewidth=.0001pt](3,-4)(3,-6)
\psline[linewidth=.0001pt](2,-5)(2,-9)
\psline[linewidth=.0001pt](1,-6)(1,-12)
\psline[linewidth=.0001pt](.5,-5)(.5,-7)
\psline[linewidth=.0001pt](0,-7)(0,-15)
\psline[linewidth=.0001pt](-.5,-6)(-.5,-10)
\psline[linewidth=.0001pt](-1.5,-7)(-1.5,-13)
\psline[linewidth=.0001pt](-2,-6)(-2,-8)
\psline[linewidth=.0001pt](-2.5,-8)(-2.5,-16)
\psline[linewidth=.0001pt](-3,-7)(-3,-11)
\psline[linewidth=.0001pt](-4,-8)(-4,-14)
\psline[linewidth=.0001pt,linestyle=dashed](-4,-4)(-4,-8)
\psline[linewidth=.0001pt,linestyle=dashed](3,0)(3,-4)
\psdot[dotsize=4pt](0,0)
\psdot[dotsize=4pt](0,2)
\psdot[dotsize=4pt](0,4)
\psdot[dotsize=4pt](0,6)
\psdot[dotsize=4pt](0,8)
\psdot[dotsize=4pt](1.5,0)
\psdot[dotsize=4pt](1.5,2)
\psdot[dotsize=4pt](1.5,4)
\psdot[dotsize=4pt](1.5,6)
\psdot[dotsize=4pt](3,0)
\psdot[dotsize=4pt](3,2)
\psdot[dotsize=4pt](3,4)
\psdot[dotsize=4pt](-1,-1)
\psdot[dotsize=4pt](-1,1)
\psdot[dotsize=4pt](-1,3)
\psdot[dotsize=4pt](-1,5)
\psdot[dotsize=4pt](-2,-2)
\psdot[dotsize=4pt](-2,0)
\psdot[dotsize=4pt](-2,2)
\psdot[dotsize=4pt](-3,-3)
\psdot[dotsize=4pt](-3,-1)
\psdot[dotsize=4pt](-4,-4)
\psdot[dotsize=4pt](-1.5,-3)
\psdot[dotsize=4pt](-.5,-2)
\psdot[dotsize=4pt](-.5,0)
\psdot[dotsize=4pt](.5,-1)
\psdot[dotsize=4pt](.5,1)
\psdot[dotsize=4pt](.5,3)
\psdot[dotsize=4pt](1,-2)
\psdot[dotsize=4pt](2,-1)
\psdot[dotsize=4pt](2,1)
\psdot[dotsize=4pt](0,-7)
\psdot[dotsize=4pt](0,-9)
\psdot[dotsize=4pt](0,-11)
\psdot[dotsize=4pt](0,-13)
\psdot[dotsize=4pt](0,-15)
\psdot[dotsize=4pt](1,-6)
\psdot[dotsize=4pt](1,-8)
\psdot[dotsize=4pt](1,-10)
\psdot[dotsize=4pt](1,-12)
\psdot[dotsize=4pt](2,-5)
\psdot[dotsize=4pt](2,-7)
\psdot[dotsize=4pt](2,-9)
\psdot[dotsize=4pt](3,-4)
\psdot[dotsize=4pt](3,-6)
\psdot[dotsize=4pt](1.5,-4)
\psdot[dotsize=4pt](.5,-5)
\psdot[dotsize=4pt](.5,-7)
\psdot[dotsize=4pt](-.5,-6)
\psdot[dotsize=4pt](-.5,-8)
\psdot[dotsize=4pt](-.5,-10)
\psdot[dotsize=4pt](-1.5,-7)
\psdot[dotsize=4pt](-1.5,-9)
\psdot[dotsize=4pt](-1.5,-11)
\psdot[dotsize=4pt](-1.5,-13)
\psdot[dotsize=4pt](-2.5,-8)
\psdot[dotsize=4pt](-2.5,-10)
\psdot[dotsize=4pt](-2.5,-12)
\psdot[dotsize=4pt](-2.5,-14)
\psdot[dotsize=4pt](-2.5,-16)
\psdot[dotsize=4pt](-4,-8)
\psdot[dotsize=4pt](-4,-10)
\psdot[dotsize=4pt](-4,-12)
\psdot[dotsize=4pt](-4,-14)
\psdot[dotsize=4pt](-3,-7)
\psdot[dotsize=4pt](-3,-9)
\psdot[dotsize=4pt](-3,-11)
\psdot[dotsize=4pt](-2,-6)
\psdot[dotsize=4pt](-2,-8)
\psdot[dotsize=4pt](-1,-5)
\uput[45](0,8){$H\cap\{x_0=4\}$}
\uput[45](1.5,6){$H\cap\{x_0=3\}$}
\uput[45](3,4){$H\cap\{x_0=2\}$}
\uput[45](3,2){$H\cap\{x_0=1\}$}
\uput[45](3,0){$H\cap\{x_0=0\}$}
\uput[45](3,-4){$H\cap\{x_0=-1\}$}
\uput[45](3,-6){$H\cap\{x_0=-2\}$}
\uput[45](2,-9){$H\cap\{x_0=-3\}$}
\uput[45](1,-12){$H\cap\{x_0=-4\}$}
\uput[45](0,-15){$H\cap\{x_0=-5\}$}
\uput[135](-4,-4){$U$}
\uput[225](-4,-8){$D$}
\end{pspicture}
\caption{: The stories of $H=\{\bullet\}$ for $S=\langle14,15,17\rangle$}
\label{fig_stories}
\end{figure}

Hence $H\cap\{x_0=i\}\neq\emptyset$ if and only if $-n_0\leq i\leq n_0-1$. We call these sets the \textbf{stories} of $H$ (see figure \ref{fig_stories}).

At first we consider the case $i\leq -1$: Then, since $\max A=c-1+g_0$, the second inequality in (\ref{neun}) always holds, hence for the underground stories,
\begin{equation}\label{zehn}H\cap\{x_0=i\}=\tilde A\cap\{x_1g_1+\ldots x_{e-1}g_{e-1}\geq -ig_0\}.\end{equation}

Now we cut $\mathds N^{e-1}$ into \textbf{strips}
\[H_0:=\{x\in\mathds N^{e-1}\vert0\leq x\cdot g<g_0\}\]
\[H_1:=\{x\in\mathds N^{e-1}\vert g_0\leq x\cdot g<2g_0\}\]
\[\vdots\]
\[H_{n_0}:=\{x\in\mathds N^{e-1}\vert n_0g_0\leq x\cdot g<(n_0+1)g_0\}\]
\[\vdots\]
of ``width'' $g_0$, and set
\[h_j:=\#\,\tilde A\cap H_j\text{, }j=0,\ldots,n_0.\]
Notice: Since $(n_0+1)g_0\geq c+g_0>\max A$, we have
\[\tilde A\cap H_{j}=\emptyset\text{ for }j\geq n_0+1.\]
Hence, by (\ref{zehn})
\[\#\,H\cap\{x_0=i\}=h_{-i}+h_{-i+1}+\ldots+h_{n_0},\]
and adding up from $-1$ to $-n_0$ yields
\[\#D=\sum_{i=-1}^{-n_0}\#\,H\cap\{x_0=i\}=0\cdot h_0+1\cdot h_1+\ldots+n_0\cdot h_{n_0}.\]

Analogously, for the stories above the ground
\[\#U=\sum_{i=0}^{n_0-1}\#\,H\cap\{x_0=i\}=n_0\cdot h_0'+(n_0-1)\cdot h_1'+\ldots+0\cdot h_{n_0}'\]
with the strips
\[H_{n_0}':=\{x\in\mathds N^{e-1}\vert c-1<x\cdot g\leq c+g_0-1\}\]
\[H_{n_0-1}':=\{x\in\mathds N^{e-1}\vert c-g_0-1<x\cdot g\leq c-1\}\]
\[\vdots\]
\[H_0':=\{x\in\mathds N^{e-1}\vert c-n_0g_0-1<x\cdot g\leq c-(n_0-1)g_0-1\}\]
and
\[h_i':=\#\,\tilde A\cap H_i'\text{, }i=0,\ldots ,n_0.\]
See \cite[Proposition 4.5]{E} for the corresponding decomposition of $L$.

Hence Wilf asked, whether
\[\frac{n_0\cdot h_0'+(n_0-1)\cdot h_1'+\ldots+0\cdot h_{n_0}'}{0\cdot h_0+1\cdot h_1+\ldots+n_0\cdot h_{n_0}}\geq\frac1{e-1}.\]
In particular, if $g_0$ divides $c$, i.\,e. $c=n_0g_0$, then the strips $H_i$ and $H_i'$ coincide and we get
\[w(S)=\frac{n_0\cdot h_0+(n_0-1)\cdot h_1+\ldots+0\cdot h_{n_0}}{0\cdot h_0+1\cdot h_1+\ldots+n_0\cdot h_{n_0}}.\]

\subsection{The Wilf ratio of an artinian $\mathds N^{e-1}$-ideal with respect to a weight vector}

Let $e\geq2$, $\gamma=(\gamma_1,\ldots,\gamma_{e-1})\in\mathds Q^{e-1}, 1\leq\gamma_1\leq\ldots\leq\gamma_{e-1}$, and $B\neq\{0\}$ the complement of an $\mathds N^{e-1}$-ideal, whose coordinate ring $R(B)$ is artinian. We define
\[n_0(\gamma,B):=\lfloor\max(B\cdot\gamma)\rfloor,\]
\[H_n(\gamma):=\{x\in\mathds N^{e-1}\vert\lfloor x\cdot\gamma\rfloor=n\}, h_n(\gamma,B):=\#\,B\cap H_n(\gamma),\]
\[H_n'(\gamma,B):=\{x\in\mathds N^{e-1}\vert n_0(\gamma,B)-\lfloor\max(B\cdot\gamma)-x\cdot\gamma\rfloor=n\}, h_n'(\gamma,B):=\#\,B\cap H_n'(\gamma,B),\]
for $n=0,\ldots,n_0(\gamma,B)$.
We call
\[w(\gamma,B):=\frac{n_0(\gamma,B)\cdot h_0'(\gamma,B)+\ldots+0\cdot h_{n_0(\gamma,B)}'(\gamma,B)}{0\cdot h_0(\gamma,B)+\ldots+n_0(\gamma,B)\cdot h_{n_0(\gamma,B)}(\gamma,B)}\]
the \textbf{Wilf ratio of the ideal} $\mathbf {\mathbf N^{e-1}\setminus B}$ \textbf{with respect to} $\mathbf \gamma$.

In the following remark we shall relate Wilf's question to our more general considerations:

\begin{remark}With the notation of section \ref{subsect31} for $B=\tilde A$ and $\gamma=g_0^{-1}g$ we have
\[n_0=n_0(\gamma,\tilde A), c+g_0-1=g_0\cdot\max(\tilde A\cdot\gamma),\]
\[H_n=H_n(\gamma), h_n=h_n(\gamma,\tilde A),\]
\[H_n'=H_n'(\gamma,\tilde A), h_n'=h_n'(\gamma,\tilde A), n=0,\ldots,n_0,\]
and
\[w(S)=w(\gamma,\tilde A).\]
Hence the question of Wilf is if
\[w(\gamma,\tilde A)\geq\frac1{e-1}.\]
\hfillplus$\square$
\end{remark}

By \cite[Theorem 20]{FGH} the Wilf ratio of $S$ is at least $\frac1{t(S)}$. In fact, if $f_1,\ldots f_t$ are the pseudo-Frobenius numbers of $S$, then $\mathds N\setminus S\subseteq\bigcup_{i=1}^t(f_i-L)$, hence $\#(\mathds N\setminus S)\leq t(S)\cdot\#L$. Analogously we will see

\begin{proposition}\label{prop_type_Wilf}

\begin{enumerate}

\item[a)]
\[w(\gamma,B)\geq\frac1{t(R(B))}.\]

\item[b)] If $B=Q_m$ is a cuboid, then $w(\gamma,B)=1$.

\end{enumerate}

\end{proposition}

\textbf{Proof} Denote by $d(B)$ and $u(B)$ resp. the denominator and the numerator of $w(\gamma,B)$ in the above formula. Then

\begin{equation}\label{dreizehn_neu}d(B)=\sum_{b\in B}\lfloor b\cdot\gamma\rfloor\text{ and }u(B)=\sum_{b\in B}\lfloor\max(B\cdot\gamma)-b\cdot\gamma\rfloor.\end{equation}

b) The bijection $Q_m\to Q_m, q\mapsto m-q$ yields $d(Q_m)=u(Q_m)$.

a) Let $m_1,\ldots,m_t$ be the maximum elements of $B$. Then $t=t(R(B))$ and $B=\bigcup_jQ_{m_j}$. Applying (\ref{dreizehn_neu}) to $B$ and to $Q_{m_j}$ for $j=1,\ldots,t$ we get
\[d(B)=d(\bigcup_jQ_{m_j})\leq\sum_{j=1}^td(Q_{m_j})=\sum_{j=1}^tu(Q_{m_j})\leq t\cdot u(B).\]\hfillplus$\square$

\begin{coro}\label{coro33ganzneu} If $e=3$ and $B$ induces a periodic tesselation, then
\[w(\gamma,B)\geq\frac1{e-1}.\]
This is immediate from \ref{coro_tess_three_equiv} and \ref{prop_type_Wilf}.\hfillplus$\square$\end{coro}

\vspace{.5cm}Coming back to Zhai's formula (\ref{Bound_Zhai}) we will see in a moment

\begin{proposition}\label{prop33} For $\gamma\in\mathds N_{\geq1}^{e-1}$, we have
\[w(\gamma,B)\geq\frac1{e-1}.\]

\end{proposition}

\textbf{Proof} Here, by (\ref{dreizehn_neu}), we have for $u=u(B)$ and $d=d(B)$
\[u+d=\max(B\cdot\gamma)\cdot \#B\text{ and }d=\sum_{b\in B}b\cdot\gamma,\]
hence by (\ref{Bound_Zhai})
\[d\leq\frac{e-1}e(u+d),\text{ equivalently }\frac ud\geq\frac1{e-1}.\]
\hfillplus$\square$

In other words, in case $\gamma\in\mathds N_{\geq1}^{e-1}$, for $i=0,\ldots,m:=\max(B\cdot\gamma)$
\[h_i:=h_i(\gamma,B)=h_i'(\gamma,B)\]
is the number of elements $b\in B$ of weight $b\cdot\gamma=i$. Hence $(h_n)_{n\in\mathds N}$ is the Hilbert function of the positively graded algebra
\[R(B)=\mathds C[X_1,\ldots,X_{e-1}]/I(\mathds N^{e-1}\setminus B), \deg X_i=\gamma_i, i=1,\ldots,e-1,\]
where
\[h_m\neq0\text{ and }h_n:=0\text{ for }n>m.\]
By \ref{prop33} Proposition,
\begin{equation}\label{elf}0\cdot h_0+1\cdot h_1+\ldots+m\cdot h_m\leq (e-1)(m\cdot h_0+(m-1)\cdot h_1+\ldots+0\cdot h_m).\end{equation}
More generally, let
\[R=\bigoplus_{i=0}^mR_i=\mathds C[X_1,\ldots,X_{e-1}]/I, R_0=\mathds C, R_m\neq0, \deg X_i=\gamma_i, i=1,\ldots,e-1\]
be an arbitrary positively graded artinian $\mathds C$-algebra with Hilbert function $h_n=\dim_{\mathds C}R_n$ and $B=\{b\in\mathds N^{e-1}\vert\underline X^b\not\in\operatorname{In}_{LEX}(I)\}$. Then by a theorem of Macaulay, see \cite[Theorem 6.1.4]{HH}, $R$ and $R(B):=\mathds C[X_1,\ldots,X_{e-1}]/I(\mathds N^{e-1}\setminus B)$ have the same Hilbert function with respect to the grading induced by $\gamma$. Hence:

\begin{remark}Formula (\ref{elf}) holds for the Hilbert function of any positively graded artinian $\mathds C$-algebra $R=\bigoplus_{i=0}^mR_i, R_m\neq0$, as well.\end{remark}

In the standard graded case $\gamma=(1,\ldots,1)$, even Eliahou's inequalities
\begin{equation}\label{zwoelf}n\cdot h_n\leq h_1(h_0+h_1+\ldots+h_{n-1}), n\in\mathds N \text{ (see \cite[Theorem 5.11]{E})}\end{equation}
hold, and $h_1\leq e-1$. Adding up from $1$ to $m$ again amounts to formula (\ref{elf}).

\vspace{.2cm}However, there are pairs $(\gamma,B)$ with $w(\gamma,B)<\frac1{e-1}$:

\begin{exampleo}

For $\gamma=(\frac32,\frac53)$, $B=\{(0,0),(1,0),(0,1),(2,0),(1,1),(0,2)\}$ one has $w(\gamma,B)=\frac5{11}$.

Here $t(R(B))=e=3$, hence there is no periodic tesselation of $\mathds Z^2$ by $B$.

\end{exampleo}

In view of \ref{tessi} and \ref{coro33ganzneu}, generalizing Wilf's question one may ask: Suppose $B$ induces a periodic tesselation of $\mathds Z^{e-1}$, do we always have
\[w(\gamma,B)\geq\frac1{e-1},\]
at least if $\#B$ is a common denominator of $\gamma_1,\ldots,\gamma_{e-1}$?

\end{document}